\documentclass[11pt,a4paper]{amsart}
\usepackage[T1]{fontenc}
\usepackage[utf8]{inputenc}
\usepackage[english]{babel}
\usepackage{lmodern}
\usepackage[margin=1in]{geometry}
\usepackage{amssymb,mathtools}
\usepackage{cite}
\usepackage{microtype}
\usepackage{enumitem}
\usepackage{booktabs}
\usepackage{tikz}
\usetikzlibrary{arrows.meta,calc,decorations.pathreplacing}
\usepackage[hidelinks,unicode]{hyperref}
\usepackage{bookmark}
\usepackage{xcolor}

\newtheorem{theorem}{Theorem}[section]
\newtheorem{lemma}[theorem]{Lemma}
\newtheorem{proposition}[theorem]{Proposition}
\newtheorem{corollary}[theorem]{Corollary}
\newtheorem{conjecture}[theorem]{Conjecture}
\theoremstyle{definition}
\newtheorem{definition}[theorem]{Definition}
\newtheorem{example}[theorem]{Example}
\newtheorem{problem}[theorem]{Problem}
\theoremstyle{remark}
\newtheorem{remark}[theorem]{Remark}
\numberwithin{equation}{section}

\newcommand{\doi}[1]{\href{https://doi.org/#1}{\texttt{doi:}\nolinkurl{#1}}}
\newcommand{\A}{\mathtt{A}}
\newcommand{\B}{\mathtt{B}}
\setlist[enumerate]{label=\textup{(\roman*)},leftmargin=2.1em,itemsep=2pt}
\newcommand{\parhead}[1]{\par\smallskip\noindent\textit{#1}\enspace}
\makeatletter
\newcommand{\jelclass}[1]{\gdef\@jelclass{#1}}
\def\@setsubjclass{{\itshape\subjclassname.}\enspace\@subjclass\@addpunct.\par
  {\itshape JEL Classification.}\enspace\@jelclass\@addpunct.}
\makeatother

\colorlet{cA}{blue!70!black}
\colorlet{cB}{red!70!black}
\colorlet{cC}{green!50!black}

\title[How long can a draft be fair?]{How long can a draft be fair?\texorpdfstring{\\}{ }Picking orders that guarantee EFX for every Borda profile}

\author[M. Anholcer]{Marcin Anholcer}
\address{Institute of Informatics and Quantitative Economics, Pozna\'n University of Economics and Business, Pozna\'n, Poland}
\email{marcin.anholcer@ue.poznan.pl}

\author[M. Bartkowiak]{Maciej Bartkowiak}
\address{Doctoral School, Pozna\'n University of Economics and Business, Pozna\'n, Poland}
\email{maciej.bartkowiak@ue.poznan.pl}

\author[B. Bosek]{Bart\l omiej Bosek}
\address{Institute of Theoretical Computer Science, Faculty of Mathematics and Computer Science, Jagiellonian University, Krak\'ow, Poland}
\email{bartlomiej.bosek@uj.edu.pl}

\author[J. Grytczuk]{Jaros\l aw Grytczuk}
\address{Faculty of Mathematics and Information Science, Warsaw University of Technology, Warsaw, Poland}
\email{jaroslaw.grytczuk@pw.edu.pl}

\thanks{This research was funded by the National Science Centre, Poland (grant number: 2023/51/B/HS4/00829).}
\subjclass[2020]{91B32, 91B14, 05C90, 68R15}
\jelclass{C78, D63, D71}
\keywords{Fair division, drafts, picking sequences, EFX, envy-freeness on graphs, Borda scores, Thue--Morse word}
\date{}

\hypersetup{pdftitle={How long can a draft be fair? Picking orders that guarantee EFX for every Borda profile},pdfauthor={Marcin Anholcer, Maciej Bartkowiak, Bartłomiej Bosek, Jarosław Grytczuk},pdfkeywords={fair division, drafts, picking sequences, EFX, Borda scores, Thue-Morse word}}

\begin{document}
\begin{abstract}
In a draft, agents take turns in a fixed order, and at each turn the current agent takes her favourite remaining item. Drafts are simple, transparent and need no preference reports. We ask which orders are fair \emph{whatever the preferences are}. Agents rank the $m$ items and value them by the Borda scores $m,m-1,\dots,1$, and envy is felt only along the edges of a graph $G$. An order is \emph{universal} if it yields an allocation that is envy-free up to any item (EFX) along every edge, for every profile.

We prove an exact formula for the worst EFX violation of a given order: it is the violation under a common ranking plus an explicit gain for each item taken by a third agent. Thus universality can be tested in linear time, although $(m!)^n$ profiles are involved. For two agents identical preferences are the worst case, and the orders that are fair after every pick are exactly the concatenations of $\A\B\B\A$ and $\B\A\A\B$, the Thue--Morse word among them. Universal orders of every length exist if and only if $G$ has an isolated vertex or a component $K_2$; otherwise every universal order has fewer than $17n+153\,\alpha(G)$ items. For complete graphs we determine all attainable lengths; for $n\ge5$ agents they are $0,1,\dots,n+1$, $2n-1$ and $2n$. In particular, the snake draft $1\,2\cdots n\,n\cdots2\,1$ known from fantasy leagues is EFX for every profile, and for $n\ge5$ no fixed order can handle more items.
\end{abstract}

\maketitle

\section{Introduction}\label{sec:intro}

Children choosing teams in a schoolyard, sports leagues distributing new players~\cite{BramsStraffin1979}, students of a business school selecting courses~\cite{BudishCantillon2012} and heirs dividing an estate~\cite{BramsTaylor1999} often use the same protocol: a \emph{draft}. The agents pick one item at a time in a prescribed order, each taking her favourite among the items still available~\cite{KohlerChandrasekaran1971}. A draft needs no money, no preference reports and no mediator~\cite{BouveretLang2011}, and everybody understands it. But the order must be fixed \emph{before} the preferences are known. Which orders are fair for every possible profile?

\parhead{EFX.} Our fairness notion is envy-freeness up to any item (EFX)~\cite{CaragiannisEtAl2019}: nobody should prefer another agent's bundle after the removal of \emph{any} single item from it. Whether EFX allocations always exist for additive valuations is one of the central open problems of fair division~\cite{Procaccia2020,AmanatidisEtAl2023}; they exist for three agents~\cite{ChaudhuryGargMehlhorn2024} and in several special cases, for instance when all agents rank the items identically~\cite{PlautRoughgarden2020} or lexicographically~\cite{HosseiniSikdarVaishXia2021}. A universal order would produce EFX allocations in the most oblivious way imaginable, without looking at the preferences at all. We show how far this obliviousness can go: in general only to a linear number of items, and for $n\ge5$ mutually comparing agents only to $2n$ items. The weaker guarantee EF1 behaves quite differently: round robin provides it for every number of items~\cite{CaragiannisEtAl2019}. The obstruction to EFX is located precisely: it comes from items taken by \emph{third} agents.

\parhead{Borda scores.} Every agent ranks the $m$ items strictly and values them at $m,m-1,\dots,1$ in this order (positive Borda scores~\cite{BaumeisterEtAl2017,DarmannKlamler2016}). This is the natural model when only rankings are available, and all agents then have the same multiset of values; to our knowledge, even in this case the existence of EFX allocations is open for four or more agents. For arbitrary additive valuations, universal orders are easy to describe and very short (Theorem~\ref{thm:interlacing}): nobody may get a third item. Borda scores are the simplest valuations for which the question has a rich answer.

\parhead{Comparison graphs.} Agents often compare themselves only with some others: colleagues in the same team, siblings, neighbours~\cite{AbebeKleinbergParkes2017,BeynierEtAl2019Local,BredereckKaczmarczykNiedermeier2022}. So we let a graph $G$ on the agents specify who compares with whom, and we require EFX along its edges (graph-EFX in the sense of Payan, Sengupta and Viswanathan~\cite{PayanSenguptaViswanathan2023}). For $G=K_n$ this is ordinary EFX.

\parhead{The question.} A picking order is a word $W$ over the agents. It is \emph{universal} for $G$ if, for every positive-Borda profile, sincere picking along $W$ produces an allocation that is EFX along every edge of $G$. The \emph{length spectrum} $\sigma(G)$ is the set of all $m\ge0$ for which a universal word of length $m$ exists. How long can a fair draft be?

\subsection*{Results at a glance}
Table~\ref{tab:results} summarises what we know about $\sigma(G)$.

\begin{table}[ht]
\centering
\small
\begin{tabular}{@{}lll@{}}
\toprule
comparison graph $G$ on $n$ vertices & length spectrum $\sigma(G)$ & reference\\
\midrule
every graph & $\supseteq\{0,\dots,n+\alpha(G)\}\cup\{2n-1,2n\}$ & Props.~\ref{prop:frontier}, \ref{prop:snake}\\
isolated vertex or component $K_2$ & all $m\ge0$, one infinite word & Cor.~\ref{cor:dichotomy}\\
all components of order $\ge3$ & $\subseteq\{0,\dots,17n+153\,\alpha(G)-1\}$ & Thm.~\ref{thm:linear}\\
$K_3$ & $\{0,\dots,11\}$ & Thm.~\ref{thm:complete-classification}\\
$K_4$ & $\{0,\dots,5\}\cup\{7,8,9\}$ & Thm.~\ref{thm:complete-classification}\\
$K_n$, $n\ge5$ & $\{0,\dots,n+1\}\cup\{2n-1,2n\}$ & Thm.~\ref{thm:complete-classification}\\
\midrule
no isolated vertices, \emph{arbitrary} additive values & $\{0,\dots,n+\alpha(G)\}$ & Thm.~\ref{thm:interlacing}\\
connected, $n\ge3$, \emph{envy-freeness} instead of EFX & $\{0\}$ & Thm.~\ref{thm:ef}\\
\bottomrule
\end{tabular}
\medskip
\caption{Lengths of universal words; $\alpha(G)$ is the independence number.}
\label{tab:results}
\end{table}

\parhead{The exact defect.} The worst EFX violation of a word over all $(m!)^n$ profiles has a closed form (Theorem~\ref{thm:exact}): it equals the violation under a \emph{common} ranking, in which the item taken at turn $t$ is worth $m+1-t$ to everybody, plus a \emph{lift gain} for each item taken by a third agent, read off from a simple lattice walk. Hence universality can be tested in linear time. Adding the formulas for the two directions of an edge gives a \emph{pair budget} (Corollary~\ref{cor:pair-budget}), which drives all our impossibility results.

\parhead{Two agents and Thue--Morse.} With two agents there are no third parties, so the worst case is identical preferences. The words that are fair after \emph{every} pick are exactly the prefixes of concatenations of the blocks $\A\B\B\A$ and $\B\A\A\B$ (Theorem~\ref{thm:abba}); the Thue--Morse word is one of them. This gives a precise sense in which the $\A\B\B\A$ pattern of tennis tie-breaks and of proposed penalty shoot-out rules~\cite{PalaciosHuerta2012} is optimal.

\begin{figure}[ht]
\centering
\begin{tikzpicture}[x=0.6cm,y=0.44cm]
  \foreach \nm/\yy/\base/\top/\exc in {3/0/{0,...,4}/{5,6}/{7,...,11},
                                       4/-1/{0,...,5}/{7,8}/{9},
                                       5/-2/{0,...,6}/{9,10}/{},
                                       6/-3/{0,...,7}/{11,12}/{},
                                       7/-4/{0,...,8}/{13,14}/{},
                                       8/-5/{0,...,9}/{15,16}/{}}{
    \node[left] at (-0.6,\yy) {$K_{\nm}$};
    \foreach \x in {0,...,17}{\draw[gray!55] (\x,\yy) circle (3.2pt);}
    \foreach \x in \base {\fill[cA] (\x,\yy) circle (3.4pt);}
    \foreach \x in \top {\fill[cB] (\x,\yy) circle (3.4pt);}
    \foreach \x in \exc {\fill[cC] (\x,\yy) circle (3.4pt);}
  }
  \foreach \x in {0,2,...,16}{\node[below] at (\x,-5.5) {\scriptsize \x};}
  \node[below] at (8.5,-6.2) {\small number of items $m$};
  \begin{scope}[shift={(18.6,-1)}]
    \fill[cA] (0,0) circle (3.4pt); \node[right] at (0.2,0) {\scriptsize $m\le n+1$};
    \fill[cB] (0,-1) circle (3.4pt); \node[right] at (0.2,-1) {\scriptsize $m\in\{2n-1,2n\}$};
    \fill[cC] (0,-2) circle (3.4pt); \node[right] at (0.2,-2) {\scriptsize sporadic};
    \draw[gray!55] (0,-3) circle (3.2pt); \node[right] at (0.2,-3) {\scriptsize none};
  \end{scope}
\end{tikzpicture}
\caption{The length spectra $\sigma(K_n)$ for $3\le n\le8$ and $m\le17$: a filled dot means that some picking order of length $m$ is EFX for every positive-Borda profile of $n$ agents.}
\label{fig:spectra}
\end{figure}
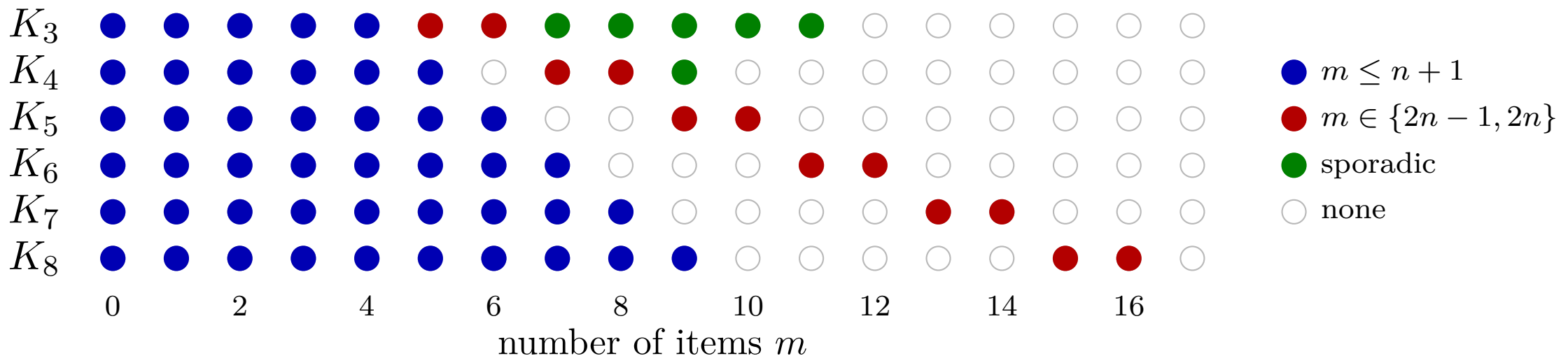

\parhead{Complete graphs.} Figure~\ref{fig:spectra} shows $\sigma(K_n)$ for small $n$. For $n\ge5$, beyond the trivial range $m\le n+1$ there are only two attainable lengths, $2n-1$ and $2n$, and the snake draft $1\,2\cdots n\,n\cdots2\,1$ attains the maximum. Moreover, for $n\ge4$ a length is attainable if and only if the block--snake draft of Chernulich and Yenmez~\cite{ChernulichYenmez2026} of that length is universal (Remark~\ref{rem:block-snake}). Repeating rounds, which is how most drafts are run, is doomed: on every graph with an edge, any concatenation of $\lceil2n/(n-2)\rceil+2$ permutation rounds fails (Corollary~\ref{cor:rounds}).

\parhead{A sincere agent is protected.} The guarantee of a universal word is individual: an agent who picks sincerely is EFX towards her neighbours whatever the others do, strategically or erratically (Remark~\ref{rem:robust}).

\parhead{Related work.} Picking sequences were introduced by Kohler and Chandrasekaran~\cite{KohlerChandrasekaran1971} and studied as elicitation-free protocols by Bouveret and Lang~\cite{BouveretLang2011}; strategic behaviour~\cite{BouveretLang2014,KalinowskiNarodytskaWalshXia2013}, welfare~\cite{AzizKalinowskiWalshXia2016}, possible and necessary outcomes~\cite{AzizWalshXia2015}, weighted agents~\cite{ChakrabortySchmidtKraepelinSuksompong2021} and the allocations reachable by some sequence~\cite{BramsKing2005,BeynierEtAl2019Sequenceability} have been investigated. Recursively balanced sequences guarantee EF1~\cite{CelineSuksompongYuen2026}, and Chernulich and Yenmez~\cite{ChernulichYenmez2026} characterise block--snake drafts by compensatory fairness axioms, including a one-sided EFX requirement towards larger bundles. Darmann and Klamler~\cite{DarmannKlamler2016} use repeated palindromes for proportional Borda allocations. Balanced alternation and the Thue--Morse word appear in work on taking turns~\cite{BramsTaylor1999,LevineStange2012,PalaciosHuerta2012,ApesteguiaPalaciosHuerta2010,AlloucheShallit1999}. What is new here is the requirement that one fixed word must give EFX for \emph{every} profile.

\section{The model and the exact defect}\label{sec:model}

\subsection{Model}
Let $G$ be a finite simple graph on the set $[n]$ of agents, and let there be $m$ goods. A \emph{positive-Borda profile} gives every agent $i$ a bijection $v_i$ from the goods onto $[m]$; values are additive. A \emph{word} $W=w_1\cdots w_m\in[n]^m$ is a picking order: at turn $t$, agent $w_t$ takes her most valuable remaining good. We name the goods by turns, so that \emph{good $t$} is the good taken at turn $t$. Agent $i$ receives $X_i=\{t:w_t=i\}$; we write $b_i=|X_i|$, $p_i(1)<\dots<p_i(b_i)$ for her turns, and $d_i=\max X_i$ for her last turn when $X_i\ne\varnothing$. The \emph{chronological value} of turn $t$ is $c(t)=m+1-t$, and $c(S)=\sum_{t\in S}c(t)$. The number $F_a=c(X_a)$ is the value of $a$'s bundle when everybody ranks the goods in the order in which they are taken.

An ordered pair of adjacent agents is an \emph{arc} $i\to j$, with \emph{observer} $i$ and \emph{target} $j$. If $X_j\ne\varnothing$, its \emph{EFX defect} in a profile is $v_i(X_j)-\min_{g\in X_j}v_i(g)-v_i(X_i)$, and its \emph{envy defect} is $v_i(X_j)-v_i(X_i)$. More generally, for an integer $h\ge0$ the \emph{EF$h$X defect} is $v_i(X_j)-v_i(E)-v_i(X_i)$, where $E$ consists of the $\min\{h,|X_j|\}$ goods of $X_j$ least valuable to $i$~\cite{AkramiRezvanSeddighin2022}; for $h=1$ and $h=0$ these are the EFX and the envy defect. Let $\delta^h_{ij}(W)$ be the largest value of the EF$h$X defect over all positive-Borda profiles, and write $\delta_{ij}=\delta^1_{ij}$. The word $W$ is \emph{universal} for $G$ if $\delta_{ij}(W)\le0$ for every arc $i\to j$ with $X_j\neq\varnothing$; comparisons with empty bundles are always fine. Agents may receive no goods. The \emph{length spectrum} $\sigma(G)$ is the set of $m\ge0$ for which some word of length $m$ is universal for $G$.

\begin{lemma}[Realizability]\label{lem:real}
A ranking $v$ of the goods occurs as the ranking of agent $i$ in a profile in which good $t$ is taken at turn $t$ if and only if $v(p)>v(g)$ for every $p\in X_i$ and every $g>p$. Such rankings are called \emph{admissible} for $i$, and the rankings of different agents can be chosen independently.
\end{lemma}
\begin{proof}
At turn $p$ all later goods were available and $i$ took good $p$; this proves necessity. Conversely, give every agent an admissible ranking. By induction on $t$, the picker at turn $t$ ranks good $t$ above all remaining goods $g>t$, so she takes it.
\end{proof}

\begin{lemma}[Upward moves]\label{lem:upward}
Fix an observer $i$. Start from the chronological ranking $c$ and move goods upwards one at a time, never moving a good above an own good of $i$ taken before it. The resulting ranking is admissible for $i$, and every own good that is never moved keeps its value $c(p)$.
\end{lemma}
\begin{proof}
A move changes the relative order of two goods only by raising the moved one. If $p\in X_i$ and $g>p$, then $p$ is initially above $g$, the good $g$ never passes $p$, and raising $p$ keeps it above $g$. If $p$ is never moved, the goods taken before $p$ stay above it, since they can only move upwards; so exactly $p-1$ goods lie above $p$.
\end{proof}

So a good may be \emph{lifted to the top} if $i$ has not picked before it was taken, and any good may be lifted to just below the last own good of $i$ taken before it.

\begin{lemma}[Chronological floor]\label{lem:floor}
The good taken at turn $p$ is worth at least $c(p)$ to the agent who takes it. Hence $v_i(X_i)\ge F_i$.
\end{lemma}
\begin{proof}
At turn $p$ only $p-1$ goods are gone, and the picker takes the best remaining one.
\end{proof}

\begin{remark}[A sincere agent is protected]\label{rem:robust}
If $W$ is universal for $G$, every agent who picks sincerely is EFX towards all her neighbours, even if the other agents take arbitrary available goods. Indeed, fix the run and the true ranking of the sincere agent. For every other agent form an auxiliary ranking that places the goods she actually takes at the top, in the order of taking. The run is sincere for the auxiliary profile, and the values of the sincere agent are unchanged; now apply universality.
\end{remark}

\subsection{The exact defect}
For $h\ge0$ let $D^h_j$ be the set of the last $\min\{h,b_j\}$ turns of $j$; thus $D^0_j=\varnothing$ and $D^1_j=\{d_j\}$. This is the set of target goods that may be deleted (none for envy-freeness, the least valuable one for EFX, which in the worst case turns out to be the last one).

\begin{definition}[Walk of an arc, lift gains]\label{def:lift}
Let $i\to j$ be an arc with $X_j\neq\varnothing$, and let $h\ge0$. Give each turn the weight
\[
\varepsilon(t)=\begin{cases}-1,&t\in X_i\quad\text{(own goods)},\\ +1,&t\in X_j\setminus D^h_j\quad\text{(retained target goods)},\\ \phantom{+}0,&\text{otherwise},\end{cases}
\]
and let $H(y)=\sum_{t\le y}\varepsilon(t)$, with $H(0)=0$, be the \emph{walk} of the arc. A turn $g\notin X_i\cup X_j$ is a \emph{third-party} turn; its \emph{lift gain} is the largest future rise of the walk,
\[
\lambda^h_{ij}(g)=\max_{g\le y\le m}\bigl(H(y)-H(g)\bigr)\ \ge0,\qquad\text{and}\qquad\Lambda^h_{ij}=\sum_{g\notin X_i\cup X_j}\lambda^h_{ij}(g).
\]
For $h=1$ we write $\lambda_{ij}$ and $\Lambda_{ij}$.
\end{definition}

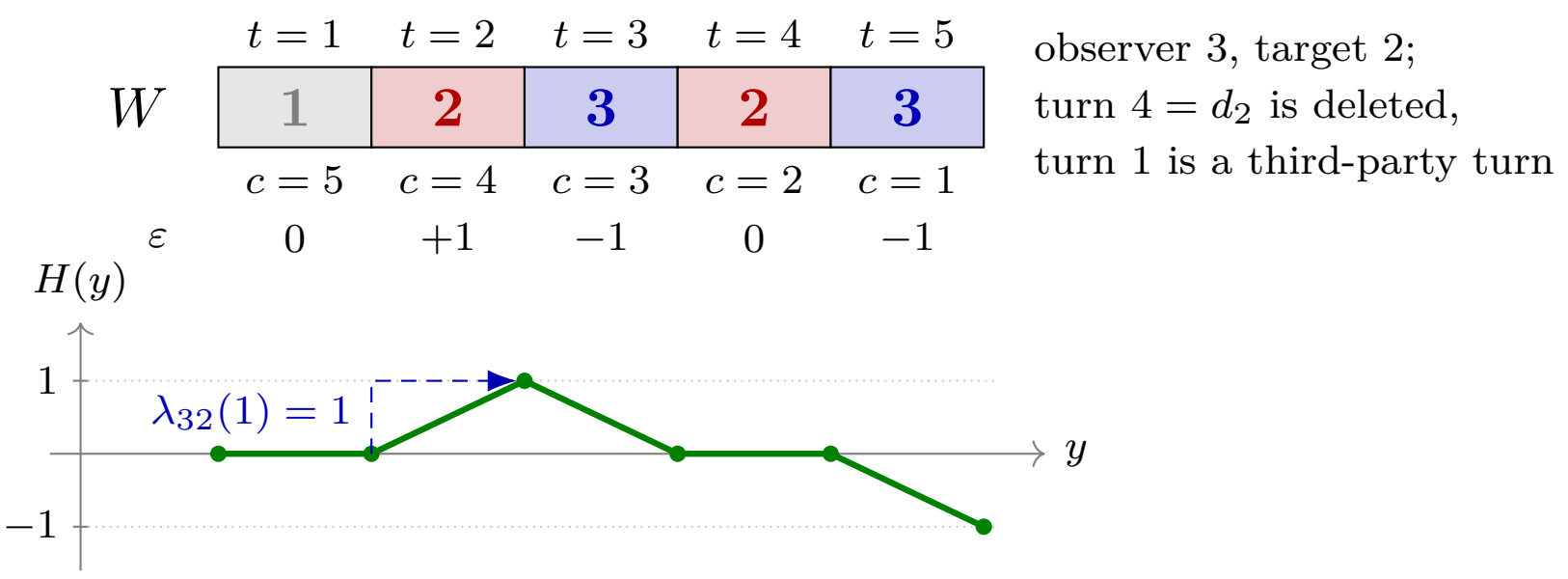
\begin{figure}[ht]
\centering
\begin{tikzpicture}[x=1.05cm,y=0.5cm]
  \foreach \t/\a/\col in {1/1/gray,2/2/cB,3/3/cA,4/2/cB,5/3/cA}{
    \fill[\col!20] (\t-0.5,0) rectangle (\t+0.5,1.1);
    \draw (\t-0.5,0) rectangle (\t+0.5,1.1);
    \node[text=\col] at (\t,0.55) {\bfseries \a};
    \node[above] at (\t,1.1) {\scriptsize $t=\t$};
    \pgfmathtruncatemacro{\cc}{6-\t}
    \node at (\t,-0.45) {\scriptsize $c=\cc$};
  }
  \node[left] at (0.3,0.55) {$W$};
  \foreach \t/\e in {1/0,2/+1,3/-1,4/{0},5/-1}{\node at (\t,-1.25) {\scriptsize $\e$};}
  \node[left] at (0.3,-1.25) {\scriptsize $\varepsilon$};
  \node[right,align=left] at (5.7,0.55) {\scriptsize observer $3$, target $2$;\\[-2pt]\scriptsize turn $4=d_2$ is deleted,\\[-2pt]\scriptsize turn $1$ is a third-party turn};
  \begin{scope}[shift={(0,-4.2)}]
    \draw[->,gray] (-0.6,0) -- (5.9,0) node[right,black] {\scriptsize $y$};
    \draw[->,gray] (-0.4,-1.6) -- (-0.4,1.8) node[above,black] {\scriptsize $H(y)$};
    \foreach \y in {-1,1}{\draw[gray] (-0.45,\y) -- (-0.35,\y) node[left,black,xshift=-2pt] {\scriptsize $\y$};
      \draw[gray!40,densely dotted] (-0.35,\y) -- (5.6,\y);}
    \draw[very thick,cC] (0.5,0) -- (1.5,0) -- (2.5,1) -- (3.5,0) -- (4.5,0) -- (5.5,-1);
    \foreach \x/\y in {0.5/0,1.5/0,2.5/1,3.5/0,4.5/0,5.5/-1}{\fill[cC] (\x,\y) circle (1.6pt);}
    \draw[densely dashed,cA,-{Latex[length=2mm]}] (1.5,0) -- (1.5,1) -- (2.45,1);
    \node[cA,left] at (1.5,0.55) {\scriptsize $\lambda_{32}(1)=1$};
  \end{scope}
\end{tikzpicture}
\caption{The walk of the arc $3\to2$ for $W=12323$. After the third-party turn $1$ it rises by at most one, so $\lambda_{32}(1)=1$ and, by~\eqref{eq:exact}, $\delta_{32}=F_2-F_3-c(d_2)+\lambda_{32}(1)=6-4-2+1=1$.}
\label{fig:walk}
\end{figure}

\begin{example}\label{ex:12323}
Take the path $1$--$2$--$3$, the word $W=1\,2\,3\,2\,3$ and the arc $3\to2$ (Figure~\ref{fig:walk}). The observer owns the turns $3,5$ and the target owns $2,4$; the last target turn $4$ is deleted, and turn $1$ is the only third-party turn. After turn $1$ the walk moves $+1,-1,0,-1$, so $\lambda_{32}(1)=1$. Here is what this means. Under the common ranking agent $3$ gets $c(3)+c(5)=4$, and the retained good $2$ is worth $c(2)=4$ to her: the defect is $0$. But agent $3$ has not picked before turn $2$, so she may rank good $2$ first and good $1$ second (Lemma~\ref{lem:upward}). This moves the third-party good $1$ below one retained target good and below no own good, and the defect becomes $5-4=1$. So $W$ is not universal for $P_3$, although it gives EFX whenever all agents have the same ranking.
\end{example}

\begin{theorem}[Exact defect]\label{thm:exact}
For every word $W$, every arc $i\to j$ with $X_j\neq\varnothing$ and every $h\ge0$,
\begin{equation}\label{eq:exact}
\delta^h_{ij}(W)=F_j-F_i-c(D^h_j)+\Lambda^h_{ij}.
\end{equation}
In particular, $W$ is universal for $G$ if and only if $F_j-F_i+\Lambda_{ij}\le c(d_j)$ for every arc $i\to j$ with $X_j\ne\varnothing$, and this can be tested in time $O(m)$ per arc.
\end{theorem}

The idea is this. A Borda value is one plus the number of goods ranked below, so the defect is a sum, over all goods $s$, of the total weight ranked above $s$. The goods of $i$ and $j$ contribute most when they are ranked in turn order, and each third-party good can be pushed down past the goods of the pair until the walk reaches its maximum.

\begin{proof}
Write $P=X_i$, $Q=X_j$, $D=D^h_j$, and let $\varepsilon,H$ be as in Definition~\ref{def:lift}; put $U(g)=\max_{y\ge g}H(y)$. By Lemma~\ref{lem:real}, $\delta^h_{ij}(W)$ is the maximum, over rankings $v$ admissible for $i$, of $v(Q)-v(E)-v(P)$, where $E$ consists of the $\min\{h,|Q|\}$ least valuable goods of $Q$.

\emph{Step 1: sorting the target.} Reassign the values of the goods of $Q$ so that they decrease in turn order. If $s$ goods of $Q$ follow an own good $p$, they now receive the $s$ smallest $Q$-values, the largest of which is at most the largest $Q$-value that followed $p$ before, and hence below $v(p)$. So admissibility is preserved, and $v(P)$ and the multiset of $Q$-values are unchanged. For such \emph{sorted} rankings $E=D$, and the objective becomes $\Phi(v)=\sum_t\varepsilon(t)v(t)$.

\emph{Step 2: a counting identity.} Since $v(t)=1+|\{s:v(s)<v(t)\}|$,
\[
\Phi(v)=\Bigl[\sum_t\varepsilon(t)+\sum_{s\in P\cup Q}A_v(s)\Bigr]+\sum_{g\notin P\cup Q}A_v(g),\qquad A_v(s)=\sum_{t:\,v(t)>v(s)}\varepsilon(t),
\]
so $A_v(s)$ is the total weight ranked above $s$. Third-party goods have weight $0$, so the bracket depends only on the relative order of $P\cup Q$.

\emph{Step 3: the bracket is at most $B=\sum_{s\in P\cup Q}H(s)$.} In turn order $A_v(s)=H(s)-\varepsilon(s)$ for $s\in P\cup Q$, and the bracket equals $B$. In general, the bracket is $\sum_t\varepsilon(t)$ plus, for every pair of goods of $P\cup Q$, the weight of the higher one. Own goods are in turn order by admissibility, target goods because $v$ is sorted, and an own good is above every later target good. So the only pairs out of turn order are a target good $q$ and a later own good $p$ ranked above it; such a pair contributes $\varepsilon(p)=-1$ instead of $\varepsilon(q)\ge0$.

\emph{Step 4: $A_v(g)\le U(g)$ for $g\notin P\cup Q$.} Let $y$ be the latest turn among $g$ and the retained target goods ranked above $g$. Every own turn $p\le y$ is above $g$: for $p<g$ by admissibility, and for $g<p\le y$ because then $y$ is a target good above $g$ and $v(p)>v(y)>v(g)$. Every retained target good above $g$ is taken by turn $y$. Hence $A_v(g)\le H(y)\le U(g)$.

\emph{Step 5: all bounds are attained.} Keep $P\cup Q$ in turn order, choose $y_g\ge g$ with $H(y_g)=U(g)$, and insert each third-party good $g$ immediately below the goods of $P\cup Q$ taken up to turn $y_g$. This ranking is sorted, and it is admissible, since a third-party good $g>p$ lies below all goods of the pair up to turn $y_g\ge g>p$. Above $g$ stand the goods of the pair up to $y_g$ and some goods of weight $0$, so $A_v(g)=U(g)$. Thus $\delta^h_{ij}(W)=B+\sum_{g\notin P\cup Q}U(g)$.

\emph{Step 6: comparison with the common ranking.} The chronological ranking $c$ is sorted, admissible, keeps $P\cup Q$ in turn order and has $A_c(g)=H(g)$. Hence $F_j-F_i-c(D)=\Phi(c)=B+\sum_gH(g)$, and subtracting, $\delta^h_{ij}(W)-(F_j-F_i-c(D))=\sum_g(U(g)-H(g))=\Lambda^h_{ij}$.
\end{proof}

\begin{remark}[Reading backwards]\label{rem:backward}
Formula~\eqref{eq:exact} is evaluated by one scan from the end of the word. Start with $d=u=0$ and read the turns $t=m,m-1,\dots,1$; the $k$-th turn read has $c(t)=k$. At an own turn put $d\leftarrow d-k$ and $u\leftarrow\max(u-1,0)$; at a retained target turn $d\leftarrow d+k$ and $u\leftarrow u+1$; at a third-party turn $d\leftarrow d+u$; the deleted target turns (the first $\min\{h,b_j\}$ turns of $j$ read) change nothing. Here $u$ is the current lift gain. At the end $d=\delta^h_{ij}(W)$. Note that the scan never uses $m$: reading backwards, all lengths are processed at once. We use this for exhaustive searches in Section~\ref{sec:complete}.
\end{remark}

\subsection{The pair budget}
\begin{corollary}[Pair budget]\label{cor:pair-budget}
Let $\{i,j\}$ be an edge with $X_i,X_j\ne\varnothing$, let $h\ge0$, and suppose $\delta^h_{ij}(W),\delta^h_{ji}(W)\le0$. Call $\kappa(g)=\lambda^h_{ij}(g)+\lambda^h_{ji}(g)$ the \emph{charge} of a third-party turn $g$. Then
\begin{equation}\label{eq:budget}
\sum_{g\notin X_i\cup X_j}\kappa(g)\le c(D^h_i)+c(D^h_j).
\end{equation}
Moreover:
\begin{enumerate}[label=\textup{(\alph*)}]
\item \emph{(next-turn charge)} if the first turn of $X_i\cup X_j$ after $g$ is not in $D^h_i\cup D^h_j$, then $\kappa(g)\ge1$;
\item \emph{(cut charge)} if $h=1$, then for every turn $s\ge g$, with $N_a=|X_a\cap(s,m]|$,
\[
\kappa(g)\ge N_i-N_j-\bigl|\{d_i\}\cap(g,m]\bigr|-\bigl|\{d_j\}\cap(g,s]\bigr|.
\]
\end{enumerate}
\end{corollary}
\begin{proof}
Adding~\eqref{eq:exact} for $i\to j$ and $j\to i$, the potentials $F_i,F_j$ cancel:
$0\ge\delta^h_{ij}+\delta^h_{ji}=-c(D^h_i)-c(D^h_j)+\sum_g\kappa(g)$. For (a), let $y$ be the first turn of $X_i\cup X_j$ after $g$. If $y$ is a retained turn of $j$, the walk of $i\to j$ rises by one at $y$, so $\lambda^h_{ij}(g)\ge1$; symmetrically if $y$ is a retained turn of $i$. For (b), evaluate $\lambda_{ji}(g)$ at the end of the word and $\lambda_{ij}(g)$ at the cut $s$:
\begin{align*}
\lambda_{ji}(g)&\ge|X_i\cap(g,m]|-|\{d_i\}\cap(g,m]|-|X_j\cap(g,m]|,\\
\lambda_{ij}(g)&\ge|X_j\cap(g,s]|-|\{d_j\}\cap(g,s]|-|X_i\cap(g,s]|.
\end{align*}
In the sum, the turns of $i$ and $j$ in $(g,s]$ cancel, and what remains is the claimed bound.
\end{proof}

In words: every third-party good pays for what happens between $i$ and $j$ after it, and the pair can pay at most the value of its last goods. If both $i$ and $j$ pick near the end, the budget is tiny. For instance, if $T$ is a turn of $i$ or $j$ after which both still pick, then at most $c(d_i)+c(d_j)$ goods taken before $T$ belong to third agents. For $h=0$ the budget is zero; this will make envy-freeness hopeless (Theorem~\ref{thm:ef}).

\begin{corollary}[Final runs]\label{cor:final-run}
Let $W$ be universal for $G$, and suppose that its last $r\ge2$ turns belong to an agent $z$ and turn $m-r$ to a neighbour $x$ of $z$. Then $(r-2)(m-b_z-b_x)\le r+2$.
\end{corollary}
\begin{proof}
All $m-b_z-b_x$ goods outside $X_z\cup X_x$ are taken before turn $m-r=d_x$. For each of them, the cut charge with $s=m-r$, $i=z$, $j=x$ is at least $r-0-1-1$, while the budget is $c(m)+c(m-r)=r+2$.
\end{proof}

\section{Two agents: ABBA and the Thue--Morse word}\label{sec:two}

Call the agents $\A$ and $\B$ and put $C=F_\B-F_\A$. With $e_t=+1$ for $\B$-turns, $e_t=-1$ for $\A$-turns and $S_k=\sum_{t\le k}e_t$, summation by parts gives $C=\sum_te_tc(t)=\sum_{k=1}^mS_k$, the signed area under the walk $S$.

\begin{theorem}[Two agents]\label{thm:two}
Let $W\in\{\A,\B\}^m$, let $i\to j$ be an arc with $X_j\ne\varnothing$, and $h\in\{0,1\}$. Then $\delta^h_{ij}(W)=F_j-F_i-c(D^h_j)$, which is the defect when both agents rank the goods identically. Consequently, $W$ is universal if and only if it yields EFX whenever both agents have the same ranking, if and only if
\begin{equation}\label{eq:two}
-c(d_\A)\le C\le c(d_\B),
\end{equation}
where an inequality involving an absent letter is omitted. The word is envy-free for every profile if and only if $C=0$.
\end{theorem}
\begin{proof}
There are no third-party goods, so this is Theorem~\ref{thm:exact} with $\Lambda=0$; under identical rankings good $t$ is worth $c(t)$ to both agents.
\end{proof}

Example~\ref{ex:12323} shows that the theorem fails for three agents.

\begin{theorem}[Fair after every pick]\label{thm:abba}
A finite or infinite word over $\{\A,\B\}$ has all its prefixes universal if and only if it is a prefix of a concatenation of the blocks $\A\B\B\A$ and $\B\A\A\B$. There are exactly $2^{\lceil m/4\rceil}$ such words of length $m$. In particular, every prefix of the Thue--Morse word $\A\B\B\A\B\A\A\B\B\A\A\B\A\B\B\A\cdots$ is universal, and at the end of each complete block the allocation is even envy-free for every profile.
\end{theorem}
\begin{proof}
\emph{If.} Inside a block starting at height $0$ the walk $S$ is $-1,0,1,0$ for $\A\B\B\A$ and $1,0,-1,0$ for $\B\A\A\B$. So complete blocks contribute $0$ to $C$, a partial block contributes $-1$, $0$ or $1$, and~\eqref{eq:two} holds because chronological values are at least~$1$. At block ends $C=0$.

\emph{Only if.} Suppose every prefix is universal and the first $4k$ letters form complete blocks, so that the walk and its area vanish at $4k$. By symmetry let letter $4k+1$ be $\A$. Each prefix considered below ends with a turn of the target $j$, so $c(d_j)=1$ and Theorem~\ref{thm:two} requires the area of the walk with sign $+1$ at the turns of $j$ to be at most~$1$. If letter $4k+2$ is $\A$, the prefix of length $4k+2$ with target $\A$ has area $1+2=3$. So it is $\B$. If letter $4k+3$ is $\A$, the prefix of length $4k+3$ with target $\A$ has area $1+0+1=2$. So it is $\B$. If letter $4k+4$ is $\B$, the prefix of length $4k+4$ with target $\B$ has area $-1+0+1+2=2$. Hence the block is $\A\B\B\A$, and induction completes the proof. Every block, complete or not, is determined by its first letter, which gives the count. The Thue--Morse word is the fixed point of $\A\mapsto\A\B\B\A$, $\B\mapsto\B\A\A\B$. Finally, with $C=0$ the chronological floors of both agents equal half of the total $m(m+1)/2$, so by Lemma~\ref{lem:floor} each agent values her own bundle at least as much as the other one.
\end{proof}

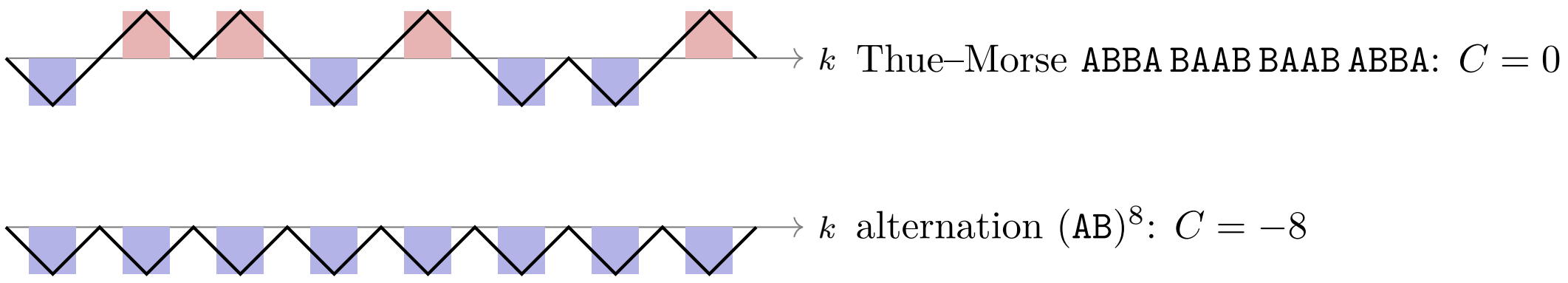
\begin{figure}[ht]
\centering
\begin{tikzpicture}[x=0.42cm,y=0.42cm]
  \def\tm{{-1,0,1,0,1,0,-1,0,1,0,-1,0,-1,0,1,0}}
  \def\rr{{-1,0,-1,0,-1,0,-1,0,-1,0,-1,0,-1,0,-1,0}}
  \foreach \nm/\arr/\yy/\lab in {tm/\tm/0/{Thue--Morse $\A\B\B\A\,\B\A\A\B\,\B\A\A\B\,\A\B\B\A$: $C=0$},rr/\rr/-3.6/{alternation $(\A\B)^8$: $C=-8$}}{
    \begin{scope}[shift={(0,\yy)}]
      \draw[gray,->] (0,0) -- (17,0) node[right,black] {\scriptsize $k$};
      \node[right] at (17.8,0) {\small \lab};
      \foreach \k in {1,...,16}{
        \pgfmathtruncatemacro{\s}{\arr[\k-1]}
        \ifnum\s>0 \fill[cB!30] (\k-0.5,0) rectangle (\k+0.5,\s);\fi
        \ifnum\s<0 \fill[cA!30] (\k-0.5,0) rectangle (\k+0.5,\s);\fi
      }
      \draw[thick] (0,0) \foreach \k in {1,...,16}{ -- (\k,{\arr[\k-1]}) };
    \end{scope}
  }
\end{tikzpicture}
\caption{The signed area under the walk $S_k$ is $C=F_\B-F_\A$. In $\A\B\B\A$ and $\B\A\A\B$ blocks it returns to $0$; under alternation it grows linearly, and~\eqref{eq:two} fails for $m=3$ and for all $m\ge5$.}
\label{fig:tm}
\end{figure}

\begin{remark}[ABBA in practice]
The Thue--Morse order was proposed for penalty shoot-outs~\cite{PalaciosHuerta2012}, where the first kicker has a measurable advantage~\cite{ApesteguiaPalaciosHuerta2010}; see also~\cite{AlloucheShallit1999,BramsTaylor1999,LevineStange2012}. Plain alternation $(\A\B)^\infty$ is universal only for $m\in\{0,1,2,4\}$.
\end{remark}

\section{Short fair drafts, the snake, and why Borda}\label{sec:short}

\begin{proposition}[Up to $n+\alpha(G)$ goods]\label{prop:frontier}
Let $m\le n+\alpha(G)$. If $m\le n$, every word in which no agent picks twice is universal for $G$. If $m>n$, choose an independent set $I$ with $|I|=m-n$, let the agents outside $I$ pick once each, and then let each agent of $I$ pick twice; this word is universal for $G$.
\end{proposition}
\begin{proof}
A comparison with a one-good bundle is always EFX. Every agent outside $I$ picks before all goods of $I$, so she values her good at least as much as each of their goods. Agents of $I$ are pairwise non-adjacent and see only singletons.
\end{proof}

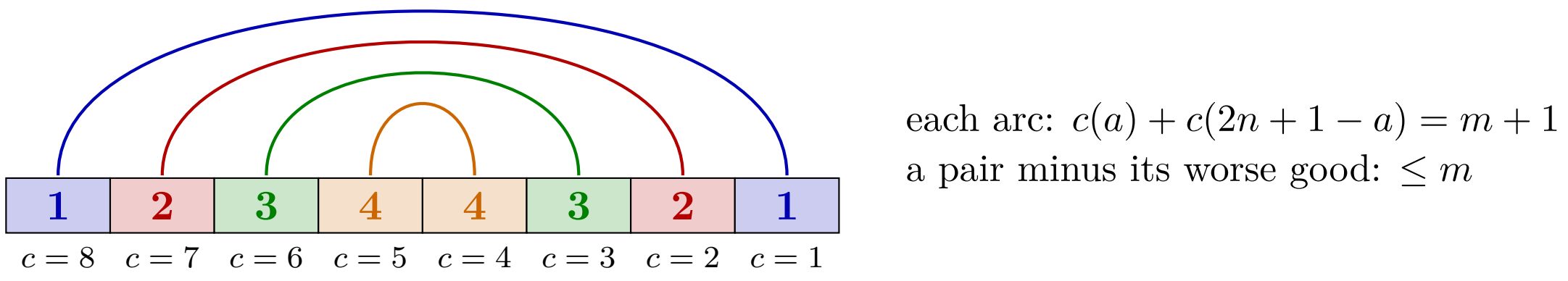
\begin{figure}[ht]
\centering
\begin{tikzpicture}[x=0.95cm,y=0.5cm]
  \def\ag{{1,2,3,4,4,3,2,1}}
  \foreach \t in {1,...,8}{
    \pgfmathtruncatemacro{\a}{\ag[\t-1]}
    \ifcase\a \or \def\col{cA}\or \def\col{cB}\or \def\col{cC}\or \def\col{orange!80!black}\fi
    \fill[\col!20] (\t-0.5,0) rectangle (\t+0.5,1);
    \draw (\t-0.5,0) rectangle (\t+0.5,1);
    \node[text=\col] at (\t,0.5) {\bfseries \a};
    \pgfmathtruncatemacro{\cc}{9-\t}
    \node at (\t,-0.45) {\scriptsize $c=\cc$};
  }
  \foreach \a/\col in {1/cA,2/cB,3/cC,4/orange!80!black}{
    \pgfmathsetmacro{\h}{1.0+0.75*(5-\a)}
    \draw[thick,\col] (\a,1.05) .. controls (\a,1.05+\h) and (9-\a,1.05+\h) .. (9-\a,1.05);
  }
  \node[right,align=left] at (9,1.6) {\small each arc: $c(a)+c(2n+1-a)=m+1$\\ \small a pair minus its worse good: $\le m$};
\end{tikzpicture}
\caption{The snake draft for $n=4$ agents and $m=8$ goods: the chronological floor of every agent exceeds what any other bundle can be worth after removing a good.}
\label{fig:snake}
\end{figure}

\begin{proposition}[The snake]\label{prop:snake}
For every $n\ge1$ and every positive-Borda profile, the snake order $1\,2\cdots n\,n\cdots2\,1$ of length $2n$ and the order $n\,1\,2\cdots(n-1)\,(n-1)\cdots2\,1$ of length $2n-1$ yield allocations that are EFX between \emph{all} pairs of agents.
\end{proposition}
\begin{proof}
In the snake, agent $a$ picks at turns $a$ and $2n+1-a$, so by Lemma~\ref{lem:floor} her pair is worth at least $c(a)+c(2n+1-a)=m+1$ to her (Figure~\ref{fig:snake}). Every other bundle is a pair, which without its worse good is worth at most $m$. In the second order agent $n$ gets her favourite good, worth $m$, while every other bundle is a pair, worth at most $m-1$ to her without its worse good; agent $a<n$ picks at turns $a+1$ and $2n-a$ and gets at least $2m-2n+1=m$.
\end{proof}

Together, Propositions~\ref{prop:frontier} and~\ref{prop:snake} show that
\[
\{0,1,\dots,n+\alpha(G)\}\cup\{2n-1,2n\}\subseteq\sigma(G)\qquad\text{for every graph }G.
\]
The snake is far from the only universal word of length $2n$.

\begin{proposition}[Many maximal drafts]\label{prop:many}
Let $W$ be a word of length $2n$ whose second half is $n\,(n-1)\cdots1$ and whose first half is a permutation of $[n]$ in which every agent $a$ occupies a position $u_a\le a+1$. Then $W$ yields EFX between all pairs of agents for every positive-Borda profile. There are $2^{n-1}$ such words, and $n!\,2^{n-1}$ words of length $2n$ that are universal for $K_n$.
\end{proposition}
\begin{proof}
Agent $a$ picks at turns $u_a$ and $2n+1-a$, so her pair is worth at least $(2n+1-u_a)+a\ge2n=m$, while any other pair without its worse good is worth at most $m$. To count the first halves, note that the last position $n$ holds agent $n$ or agent $n-1$; in the second case, agent $n$ takes over the role of agent $n-1$ in a smaller instance. Hence the number of first halves doubles with each new agent. Relabelling the agents gives the last claim.
\end{proof}

\begin{conjecture}\label{conj:2n}
For $n\ge3$, the universal words of length $2n$ for $K_n$ are exactly the $n!\,2^{n-1}$ words of Proposition~\ref{prop:many} with relabelled agents, and the universal words of length $2n-1$ are exactly the $n!$ relabellings of $n\,1\,2\cdots(n-1)\,(n-1)\cdots2\,1$.
\end{conjecture}

We have verified Conjecture~\ref{conj:2n} by computer for $3\le n\le6$.

\subsection{Why Borda scores?}
For arbitrary additive valuations the question becomes simple, and the answer is poor. In this subsection values are arbitrary nonnegative numbers, ties are broken arbitrarily, and universality is required for every tie-breaking.

\begin{theorem}[Prefix interlacing]\label{thm:interlacing}
A word is universal for $G$ under arbitrary nonnegative additive valuations if and only if, for every arc $i\to j$ and every $r\le b_j-1$, agent $i$ picks at least $r$ times and $p_i(r)<p_j(r)$. Consequently, if $G$ has no isolated vertices, such words exist exactly for $m\le n+\alpha(G)$.
\end{theorem}
\begin{proof}
\emph{Sufficiency.} Let $t_1,\dots,t_b$ be the goods of $j$ and $o_1,\dots,o_{b-1}$ the first $b-1$ goods of $i$. Remove $t_s$. For $r<s$, the good $t_r$ was still available when $i$ took $o_r$, and for $r\ge s$ so was $t_{r+1}$; hence $v_i(o_r)\ge v_i(t_r)$, respectively $v_i(o_r)\ge v_i(t_{r+1})$. Summing gives $v_i(X_i)\ge v_i(X_j\setminus\{t_s\})$.

\emph{Necessity.} Suppose that $i$ picks fewer than $r$ times before turn $K=p_j(r)$, where $r\le b_j-1$. Let every other agent value good $t$ at $m+1-t$, and let $i$ value the goods $1,\dots,K$ at $1$ and the others at $0$, breaking ties towards smaller labels. Then good $t$ is taken at turn $t$, $v_i(X_i)\le r-1$, and $X_j$ without a good of value $0$ is worth $r$ to $i$.

\emph{Lengths.} If two neighbours both picked twice, each would have to pick first; and if $j$ picked three times, all neighbours of $j$ would have to pick twice. Hence agents picking twice form an independent set and nobody picks three times, so $m\le n+\alpha(G)$. The words of Proposition~\ref{prop:frontier} satisfy the condition.
\end{proof}

So everything beyond $n+\alpha(G)$ is due to the Borda structure: for $K_n$ the maximum grows from $n+1$ to $2n$, and on sparse graphs even more may be possible (Section~\ref{sec:open}).

\subsection{Envy-freeness is hopeless}
\begin{theorem}[Universal envy-freeness]\label{thm:ef}
Let $m>0$ and let $G$ have no isolated vertex. A word of length $m$ that is envy-free for every positive-Borda profile exists if and only if $G$ has a component $K_2$ and $m\equiv0$ or $3\pmod4$. In particular, a connected graph with at least three vertices admits no such word.
\end{theorem}
\begin{proof}
Let $W$ be such a word. Under the common chronological ranking both directions of an edge force $F_i=F_j$, so a component that receives a good gives a good to each of its vertices. For an edge $\{i,j\}$ with nonempty bundles the budget~\eqref{eq:budget} vanishes for $h=0$, so by the next-turn charge no turn of $i$ or $j$ follows a good of a third agent. If a component receiving goods contained a path $x,y,z$, the edge $xy$ would put all goods of $z$ after those of $x$, and the edge $yz$ the other way round. The same argument applied to edges of two different components shows that only one component receives goods. Hence it is a $K_2$, and $F_\A=F_\B$ requires $m(m+1)/2$ to be even. Conversely, $(\A\B\B\A)^q$ and $(\A\B\B\A)^q\A\B\B$ have equal chronological floors, so each agent gets at least half of the total value.
\end{proof}

\section{How long can a fair draft be?}\label{sec:linear}

Throughout this section every component of $G$ has at least three vertices, and $W$ is universal for $G$. Let
\[
w=\min_{\{i,j\}\in E(G)}(m-b_i-b_j)
\]
be the number of goods outside the heaviest edge. We play two estimates of $w$ against each other. On the one hand, $w$ is large: under the common ranking neighbours have nearly equal bundle values, so a vertex next to a heavy edge must be nearly as rich as its ends (Lemma~\ref{lem:outside}). On the other hand, $w$ is small: in the last $\approx w/3$ turns, two neighbours cannot both pick twice, and an agent who picks often there is expensive for the pair budget of her neighbours (Lemma~\ref{lem:window}).

\begin{lemma}[Local terminal constraints]\label{lem:local}
Let $\{i,j\}$ be an edge with $X_i,X_j\ne\varnothing$, let $1\le\ell\le m$, and let $N_i,N_j$ be the numbers of picks of $i$ and $j$ among the last $\ell$ turns.
\begin{enumerate}[label=\textup{(\alph*)}]
\item If $3(\ell-1)\le m-b_i-b_j$, then $N_i\le1$ or $N_j\le1$.
\item If $N_j\le1$ and $N_i>2$, then $(N_i-2)(m-\ell-b_i-b_j+N_i)\le m+\ell+1-N_i$.
\end{enumerate}
\end{lemma}
\begin{proof}
Put $s=m-\ell$. The number of third-party goods taken by turn $s$ is $\rho=m-\ell-b_i-b_j+N_i+N_j$.

(a) Suppose $N_i,N_j\ge2$. Then $\rho\ge3(\ell-1)-\ell+4=2\ell+1$. For each of these goods, the first turn of $i$ or $j$ after it is not a last turn, because its owner picks at least twice in the final interval. By the next-turn charge, the charges add up to at least $2\ell+1$, while both last turns lie in the final interval, so the budget is less than $2\ell$.

(b) Since $N_j\le1$, at most $1-N_j$ last turns of $j$ lie in $(g,s]$. The cut charge at $s$ charges each of the $\rho$ goods at least $N_i-N_j-1-(1-N_j)=N_i-2$. The last turn of $i$ is preceded by $N_i-1$ of her turns in the final interval, so $c(d_i)\le\ell+1-N_i$, while $c(d_j)\le m$. Compare with~\eqref{eq:budget} and use $\rho\ge m-\ell-b_i-b_j+N_i$.
\end{proof}

\begin{lemma}[Counting a final interval]\label{lem:window}
Let $m>0$ and let $t\ge1$ be an integer. If $m<\frac{2t-1}{3}w+t(t+2)$, then $w<3\bigl(n+t\,\alpha(G)\bigr)$.
\end{lemma}
\begin{proof}
Take the last $\ell=\lfloor w/3\rfloor+1$ turns and let $N_x$ count the picks of $x$ there. Call $x$ \emph{busy} if $N_x\ge2$. By Lemma~\ref{lem:local}(a), busy agents form an independent set. We claim that $N_x\le t+1$ for every $x$. Then $w/3<\ell=\sum_xN_x\le n+t\,\alpha(G)$, as required.

Suppose $N_x\ge t+2$. Since the component of $x$ has at least three vertices, $x$ has a neighbour $y$; she is not busy, and $X_y\ne\varnothing$, since otherwise $y$ would face a bundle with at least two goods. Apply Lemma~\ref{lem:local}(b) to $i=x$, $j=y$. Its left side increases with $N_x$ and its right side decreases, and $m-b_x-b_y\ge w$; so
$t(w-\ell+t+2)\le m+\ell-t-1$. Inserting $\ell\le w/3+1$ gives $t(\tfrac23w+t+1)\le m+\tfrac13w-t$, that is, $m\ge\frac{2t-1}{3}w+t(t+2)$, a contradiction.
\end{proof}

\begin{lemma}[Goods outside an edge]\label{lem:outside}
For every edge $\{i,j\}$, the number $q=b_i+b_j$ satisfies $q(q+1)\le\frac23m(m+1)+2m$. In particular, $w>3m/17$ whenever $m\ge160$.
\end{lemma}
\begin{proof}
By Theorem~\ref{thm:exact} and $\Lambda\ge0$, every arc $a\to b$ has $F_b-F_a\le c(d_b)\le m$ (trivially if $X_b=\varnothing$). Put $X=(F_i+F_j)/2$; both ends of the edge have $F\ge X-m/2$, and a third vertex $k$ adjacent to one of them has $F_k\ge X-3m/2$. The values $F$ of all agents add up to $m(m+1)/2$, so $m(m+1)/2\ge3X-3m/2$. Since the two ends hold $q$ distinct turns, $2X\ge1+2+\dots+q=q(q+1)/2$, and the inequality follows. If $q\ge14m/17$, dividing by $m$ gives $\frac{10}{867}m\le\frac{94}{51}$, that is, $m\le159.8$.
\end{proof}

\begin{theorem}[A linear length bound]\label{thm:linear}
If every component of $G$ has at least three vertices, every universal word for $G$ has
\[
m<17n+153\,\alpha(G)\le170n-153 .
\]
\end{theorem}
\begin{proof}
As $n\ge3$ and $\alpha(G)\ge1$, we may assume $m\ge204$. Lemma~\ref{lem:outside} gives $w>3m/17$, so $\frac{17}{3}w>m$, and Lemma~\ref{lem:window} applies with $t=9$. Hence $3m/17<w<3(n+9\alpha(G))$. Finally $\alpha(G)\le n-1$.
\end{proof}

\begin{corollary}[Dichotomy]\label{cor:dichotomy}
For every graph $G$ the following are equivalent:
\begin{enumerate}
\item $\sigma(G)$ contains every $m\ge0$;
\item $\sigma(G)$ is infinite;
\item some infinite word has all its prefixes universal for $G$;
\item $G$ has an isolated vertex or a component isomorphic to $K_2$.
\end{enumerate}
\end{corollary}
\begin{proof}
If (iv) holds, give all goods to the isolated vertex, or assign all turns to the $K_2$ component along the Thue--Morse word (Theorem~\ref{thm:abba}). This gives (iii) and (i), each of which implies (ii). If (iv) fails, Theorem~\ref{thm:linear} shows that (ii) fails.
\end{proof}

\begin{remark}[Constants]\label{rem:constants}
We made no effort to optimise the constants; a finer analysis, which compares bundle sizes along three-vertex paths, improves them considerably. The method also extends to EF$h$X, where $h$ goods may be removed~\cite{AkramiRezvanSeddighin2022}.
\end{remark}

\section{Complete graphs}\label{sec:complete}

On $K_n$ with $m>n$ all bundles are nonempty, since an empty agent would face a bundle with at least two goods. We first record two adversarial tests that work on every graph; both come from lifting a block of early goods to the top.

\begin{lemma}[Top block]\label{lem:top-block}
Let $i\to j$ be an arc and $0\le a\le b_i$. Let $T$ be a set of $\tau$ goods of $j$ with $b_j\ge\tau+1$, each taken before the turn $p_i(a+1)$ if $a<b_i$. Then
\[
\delta_{ij}(W)\ge(\tau-a)(m+1)-a\tau-\tbinom{\tau+1}2+\tbinom{a+1}2-\sum_{\ell=a+1}^{b_i}c\bigl(p_i(\ell)\bigr).
\]
\end{lemma}
\begin{proof}
Let $A$ be the set of the first $a$ own goods. Raise the goods of $A\cup T$, one at a time in chronological order, to the top, each just below those already raised. Every own good taken before a good of $T$ lies in $A$ and has been raised before it, so the ranking is admissible by Lemma~\ref{lem:upward}, and the own goods outside $A$ keep their chronological values. The raised goods occupy the top $a+\tau$ places in chronological order; if $A$ stood above $T$, the goods of $A$ would receive $a(m+1)-\binom{a+1}2$ and those of $T$ would receive $\tau(m+1)-a\tau-\binom{\tau+1}2$, and every good of $T$ preceding a good of $A$ only moves value from $A$ to $T$. Since $b_j>\tau$, removing the least valuable good of $X_j$ leaves at least the value of $T$.
\end{proof}

\begin{corollary}[Count gap and late-pick witness]\label{cor:tests}
Let $i\to j$ be an arc.
\begin{enumerate}[label=\textup{(\alph*)}]
\item If $b_i=p$ and $b_j\ge p+d+1$ with $d\ge0$, then $\delta_{ij}(W)\ge dm-p^2-2pd-\binom d2$. In particular, $b_j\ge p+2$ forces $m\le p(p+2)$.
\item If $1\le r\le b_i$, $b_j\ge r+1$ and $p_j(r)<p_i(r)$, then $\delta_{ij}(W)\ge\sum_{\ell=r}^{b_i}p_i(\ell)-(b_i-r)(m+1)-r^2$.
\end{enumerate}
\end{corollary}
\begin{proof}
For (a) apply Lemma~\ref{lem:top-block} with $a=p$ and $T$ the first $p+d$ goods of $j$; for (b) with $a=r-1$ and $T$ the first $r$ goods of $j$. The rest is arithmetic.
\end{proof}

For $r=1$ the witness says: if a neighbour with at least two goods picks before the observer's first turn, then a single-good observer must pick at turn~$1$, and a two-good observer must have $p_i(1)+p_i(2)\le m+2$.

\begin{theorem}[Forbidden lengths up to $3n$]\label{thm:upto3n}
No word of length $m$ is universal for $K_n$ if
\begin{enumerate}[label=\textup{(\alph*)}]
\item $n\ge4$ and $n+2\le m\le2n-2$, or
\item $n\ge4$ and $2n+2\le m\le3n-1$, or
\item $n\ge5$ and $m\in\{2n+1,3n\}$.
\end{enumerate}
\end{theorem}
\begin{proof}
(a) Let $k=m-n\ge2$. Some bundle is a singleton, so by Corollary~\ref{cor:tests}(a) all sizes are $1$ or $2$: there are $n-k$ singletons and $k$ pairs. By the case $r=1$ of the witness, all singletons pick before every pair owner, i.e.\ at turns $1,\dots,n-k$. The owner $i_0$ of turn $n-k+1$ owns at most one of the last two turns, so another pair owner $i$ has turns summing to at least $(n-k+2)+(m-1)=2n+1>m+2$; she contradicts the witness towards $i_0$.

(b), (c) Now $2n+1\le m\le3n$ and $m\ge9$. A singleton would force all sizes to be at most $2$, and a bundle of size at least $4$ would exclude pairs; so there are $k=m-2n$ triples and $n-k$ pairs. The latest second pick, at turn $t_2\ge2n$, belongs to a triple owner $T$, since otherwise the witness with $r=2$ towards a triple owner gives defect at least $t_2-4$. Let $t_1<t_2<t_3$ be the turns of $T$. If $k\ge2$, the witness with $r=2$ from $T$ towards another triple owner requires $t_2+t_3\le m+5$, but $t_2+t_3\ge4n+1>m+5$ in all claimed cases. If $k=1$, then $m=2n+1$, $n\ge5$, and $T$ owns the last two turns. Her chronological value is $m+4-t_1$, and the pair containing good $1$ is worth $m$ to her without its worse good, so $t_1\le4$. The agent who starts last picks first at a turn $q_1\ge n>4$, so she owns a pair, with second turn $q_2<t_2$. She may lift good $t_1$ to the top and good $t_2$ to just below $q_2$ (Lemma~\ref{lem:upward}); then the two better goods of $T$ are worth $m+(m-q_2)$ to her, while her own pair is worth $2m+2-q_1-q_2$, which is smaller by $q_1-2>0$.
\end{proof}

\begin{lemma}[Count spread]\label{lem:spread}
If $W$ is universal for $K_n$, $n\ge3$, and all bundles are nonempty, then the smallest size $p$ and the largest size $B$ satisfy $B\le p+\lceil p/(n-2)\rceil$.
\end{lemma}
\begin{proof}
Let $d=\lceil p/(n-2)\rceil$ and suppose $B\ge p+d+1$. Then $m\ge np+d+1$, and Corollary~\ref{cor:tests}(a) from a $p$-good owner to a $B$-good owner gives defect at least $d(np+d+1)-p^2-2pd-\binom d2=(n-2)pd-p^2+d(d+3)/2>0$.
\end{proof}

\begin{figure}[ht]
\centering
\begin{tikzpicture}[x=0.72cm,y=0.62cm]
  \def\word{{9,9,0,9,1,9,9,0,2,0,3,4,1,1}}
  \foreach \t in {1,...,14}{
    \pgfmathtruncatemacro{\a}{\word[\t-1]}
    \ifnum\a=9 \def\col{gray}\def\lab{$\ast$}\fi
    \ifnum\a=0 \def\col{cA}\def\lab{$y$}\fi
    \ifnum\a=1 \def\col{cB}\def\lab{$z$}\fi
    \ifnum\a=2 \def\col{gray}\def\lab{$x_1$}\fi
    \ifnum\a=3 \def\col{gray}\def\lab{$x_2$}\fi
    \ifnum\a=4 \def\col{gray}\def\lab{$x_3$}\fi
    \fill[\col!20] (\t-0.5,0) rectangle (\t+0.5,1);
    \draw (\t-0.5,0) rectangle (\t+0.5,1);
    \node[text=\col] at (\t,0.5) {\small \lab};
  }
  \node[left] at (0.3,0.5) {$\cdots$};
  \node[below] at (8,0) {\scriptsize $T$};
  \draw[decorate,decoration={brace,amplitude=4pt}] (8.5,1.15) -- (12.5,1.15) node[midway,above=5pt] {\scriptsize block $D$ ($L=4$)};
  \draw[decorate,decoration={brace,amplitude=4pt}] (12.6,1.15) -- (14.5,1.15) node[midway,above=5pt,xshift=4pt] {\scriptsize $z^r$ ($r=2$)};
  \draw[decorate,decoration={brace,mirror,amplitude=4pt}] (0.5,-0.15) -- (7.5,-0.15) node[midway,below=5pt] {\scriptsize each $\ast$ before $T$: charge $\ge1$};
  \draw[decorate,decoration={brace,mirror,amplitude=4pt}] (10.5,-0.15) -- (12.5,-0.15) node[midway,below=5pt,align=center] {\scriptsize $a=2$ goods:\\[-2pt]\scriptsize charge $\ge r-1$};
  \node[right,align=left] at (14.7,0.5) {\scriptsize budget of $\{y,z\}$:\\[-1pt]\scriptsize $c(d_y)+c(d_z)$\\[-1pt]\scriptsize $=(a+r+1)+1$};
\end{tikzpicture}
\caption{The end of a word in the proof of Theorem~\ref{thm:deficit}. The last goods of $y$ and $z$ are nearly worthless, yet the pair budget must pay for every good of a third agent taken before $T$; so few goods lie outside $X_y\cup X_z$.}
\label{fig:deficit}
\end{figure}
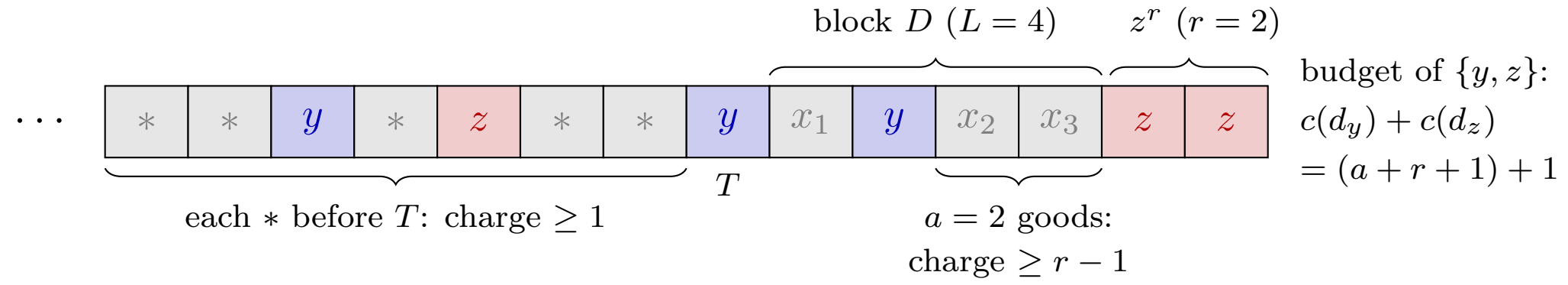

The next theorem is the heart of the classification. On a complete graph every two agents are adjacent, so the pair budget may be applied to the last picker together with an agent who picks shortly before her final run; their last goods are nearly worthless, yet their budget must pay for the goods of all other agents (Figure~\ref{fig:deficit}).

\begin{theorem}[Terminal deficit]\label{thm:deficit}
Let $n\ge3$, and let $W$ be universal for $K_n$ with $m>n$. If $B\ge C$ are the two largest bundle sizes, then
\[
m-B-C\le2n-1 .
\]
In words, the $n-2$ smallest bundles hold at most $2n-1$ goods together.
\end{theorem}
\begin{proof}
Put $w=m-B-C$, so that $w\le m-b_y-b_z$ for all $y\ne z$, and let the word end with a maximal run $z^r$.

\emph{Long final runs.} If $r\ge3$, Corollary~\ref{cor:final-run} applied to $z$ and the owner $x$ of turn $m-r$ gives $w\le(r+2)/(r-2)\le5\le2n-1$.

\emph{Short final runs.} Let $r\le2$, and let $D$ be the longest block of distinct agents other than $z$ immediately before $z^r$; its length $L$ satisfies $1\le L\le n-1$. If $D$ starts the word, then $w\le L-1$. Otherwise, by maximality, the turn $T$ just before $D$ belongs to $z$ or to a member of $D$. In the latter case call that member $y$; otherwise let $y$ be any member of $D$. Both $y$ and $z$ pick after $T$. Let $a$ be the number of turns of $D$ after the turn of $y$, and $R$ the number of goods outside $X_y\cup X_z$ taken before $T$. Each of these $R$ goods has charge at least $1$ (next-turn charge). Each of the $a$ goods of $D$ after the turn of $y$ has charge at least $\lambda_{yz}\ge r-1$, since after it the pair has only the run $z^r$ whose last turn is deleted. The budget is $c(d_y)+c(d_z)=(a+r+1)+1$. Hence $R+(r-1)a\le a+r+2$, and
\[
w\le m-b_y-b_z=R+L-1\le L+r+1+(2-r)a .
\]
For $r=1$ this is at most $L+2+(L-1)\le2n-1$, and for $r=2$ it is $L+3\le2n-1$.
\end{proof}

\begin{corollary}\label{cor:deficit}
Let $W$ be universal for $K_n$ with $m>n\ge3$, and let $p$ be the smallest bundle size.
\begin{enumerate}[label=\textup{(\alph*)}]
\item If $n\ge6$, every bundle has at most three goods, so $m\le3n$.
\item If $n=5$, then $m\le17$; if $n=4$, then $m\le15$; if $n=3$, then $m\le21$.
\end{enumerate}
\end{corollary}
\begin{proof}
By Theorem~\ref{thm:deficit}, $(n-2)p\le2n-1$, and by Lemma~\ref{lem:spread} all sizes are at most $p+\lceil p/(n-2)\rceil$.
(a) For $n\ge6$ we get $p\le2$, so all sizes are at most $3$.
(b) For $n=5$, $p\le3$, all sizes are at most $4$, and $m\le9+8$. For $n=4$, $p\le3$. If $p\le2$, all sizes are at most $3$ and not all equal $3$, so $m\le11$; if $p=3$, either some size is $5$ and Corollary~\ref{cor:tests}(a) gives $m\le15$, or $m\le7+8$. For $n=3$, $p\le5$ and $B\le2p$. If $p\le4$ then $m\le5p\le20$. If $p=5$ and $m\ge22$, then $B\ge9$, and Corollary~\ref{cor:tests}(a) with $d=B-6\in\{3,4\}$ gives defect at least $dm-25-10d-\binom d2>0$.
\end{proof}

\begin{theorem}[Complete graphs]\label{thm:complete-classification}
The length spectrum of $K_n$ is
\[
\sigma(K_n)=\begin{cases}
\mathbb Z_{\ge0}, & n\le2,\\
\{0,1,\dots,11\}, & n=3,\\
\{0,1,\dots,5\}\cup\{7,8,9\}, & n=4,\\
\{0,1,\dots,n+1\}\cup\{2n-1,2n\}, & n\ge5.
\end{cases}
\]
\end{theorem}
\begin{proof}
For $n\le2$ use Corollary~\ref{cor:dichotomy}. For every $n$, the lengths $m\le n+1$ and $m\in\{2n-1,2n\}$ are attainable by Propositions~\ref{prop:frontier} and~\ref{prop:snake}. For $n\ge6$, Theorem~\ref{thm:upto3n} excludes $n+2\le m\le2n-2$ and $2n+1\le m\le3n$, and Corollary~\ref{cor:deficit}(a) excludes $m>3n$. This settles $n\ge6$.

For $n\in\{3,4\}$ the remaining attainable lengths are witnessed by the words of Table~\ref{tab:small}, which are universal by a direct evaluation of~\eqref{eq:exact}. Theorem~\ref{thm:upto3n} excludes $m\in\{6,10,11\}$ for $n=4$ and $m\in\{7,8\}\cup\{11,\dots,15\}$ for $n=5$, and Corollary~\ref{cor:deficit}(b) leaves only finitely many lengths: $12\le m\le21$ for $n=3$, $12\le m\le15$ for $n=4$, and $m\in\{16,17\}$ for $n=5$. These are excluded by the exhaustive search of Remark~\ref{rem:search}.
\end{proof}

\begin{table}[ht]
\centering
\small
\begin{tabular}{@{}cl|rrrrrr@{}}
\toprule
$m$ & word for $K_3$ & $\delta_{12}$ & $\delta_{13}$ & $\delta_{21}$ & $\delta_{23}$ & $\delta_{31}$ & $\delta_{32}$\\
\midrule
7 & $\mathtt{1231233}$ & $-5$ & $-2$ & $-2$ & $-2$ & $-1$ & $-1$\\
8 & $\mathtt{11223333}$ & $-9$ & $0$ & $-3$ & $0$ & $-2$ & $-2$\\
9 & $\mathtt{123231312}$ & $0$ & $0$ & $-1$ & $-3$ & $-2$ & $0$\\
10 & $\mathtt{1233212311}$ & $-1$ & $-1$ & $-1$ & $-3$ & $0$ & $-2$\\
11 & $\mathtt{12311223333}$ & $-10$ & $0$ & $-1$ & $0$ & $0$ & $0$\\
\midrule
9 & \multicolumn{7}{l}{$\mathtt{123412344}$ and $\mathtt{123432144}$ for $K_4$}\\
\bottomrule
\end{tabular}
\medskip
\caption{Universal words for $K_3$ at lengths $7$--$11$, with their six defects~\eqref{eq:exact}, and for $K_4$ at length $9$.}
\label{tab:small}
\end{table}

For instance, in the last word for $K_3$ we have $F_1=26$, $F_3=19$ and $c(d_3)=1$, and the three third-party turns $2,6,7$ of agent $2$ have lift gains $2,3,3$ for the arc $1\to3$; so $\delta_{13}=19-26-1+8=0$. For $K_4$ and $m=9$, agent $4$ picks at turns $4,8,9$, so her floor is $9$, and each other agent picks once before turn $4$ and once before turn $8$, so the two goods of each of them dominate any two goods of agent $4$.

\begin{remark}[Exhaustive search]\label{rem:search}
By Remark~\ref{rem:backward}, a word can be built from its end, letter by letter, and all lengths are processed in one sweep. For every arc we keep the pair $(d,u)$ of the backward scan. The future of the scan depends only on the $u$-values and on which agents have already appeared once or at least twice (the first letter of an agent read backwards is her deleted last turn), so for each such state it suffices to keep the coordinatewise minimal vectors of $d$-values. Moreover, once $i$ and $j$ have both appeared, the sum of the $d$-values of $i\to j$ and $j\to i$ never decreases: a letter $i$ adds $-k$ to the first and $+k$ to the second, and a third-party letter adds nonnegative gains. This is the pair budget read backwards, and the same holds along every directed cycle; states violating it are discarded. Using the symmetry of $K_n$ to fix the order in which agents first appear, the search for $K_3$, $K_4$ and $K_5$ up to the lengths $21$, $15$ and $17$ takes about a second in Python and finds universal words exactly at the lengths stated in Theorem~\ref{thm:complete-classification}.
\end{remark}

\begin{remark}[Block--snake drafts]\label{rem:block-snake}
Write $m=kn+r$ with $0\le r<n$, split the agents into blocks $\mathcal A=(1,\dots,n-r)$ and $\mathcal B=(n-r+1,\dots,n)$, and let $\mathcal A_t,\mathcal B_t$ be these blocks for odd $t$ and their reversals for even $t$. The \emph{block--snake draft} of Chernulich and Yenmez~\cite{ChernulichYenmez2026} is $\mathcal A_1\mathcal B_1\mathcal A_2\mathcal B_2\cdots\mathcal A_k\mathcal B_k\mathcal B_{k+1}$. For $n\ge4$ and every $m$, the block--snake draft of length $m$ is universal for $K_n$ if and only if $m\in\sigma(K_n)$: the words of Propositions~\ref{prop:frontier} and~\ref{prop:snake} and the second word of Table~\ref{tab:small} are block--snake drafts. For $n=3$ the block--snake draft is universal for $m\le9$ but not for $m=10,11$: for $m=10$ it is $1232131233$, and agent~$3$ with the admissible values $10,9,8,6,7,5,4,3,2,1$ for the goods in turn order gets $16$, while the bundle of agent~$1$ without its worst good is worth $17$ to her. So the drafts used in practice are, for $n\ge4$, as good as any fixed order.
\end{remark}

\section{Repeating rounds does not work}\label{sec:rounds}

Most drafts are run in rounds: in each round, every agent picks the same number of times. The pair budget shows that this cannot go on for long, whatever the order within the rounds, the initial segment and the comparison graph.

\begin{theorem}[Balanced blocks]\label{thm:blocks}
Let $n\ge3$, let $\{i,j\}$ be an edge of $G$, and let $W=W_0B_1\cdots B_T$ with $T\ge2$, where $W_0$ is arbitrary and every agent occurs exactly $a_t\ge1$ times in $B_t$. If $(n-2)\sum_{t=1}^{T-2}a_t\ge2na_T$, then $W$ is not universal for $G$.
\end{theorem}
\begin{proof}
The last turns of $i$ and $j$ are distinct turns of $B_T$, so the budget is at most $2na_T-1$. A good of a third agent in a block $B_t$ with $t\le T-2$ is followed, within $B_t$ or $B_{t+1}$, by a turn of $i$ or $j$ that is not a last turn, so its charge is at least $1$. There are $(n-2)\sum_{t\le T-2}a_t$ such goods, contradicting~\eqref{eq:budget}.
\end{proof}

\begin{corollary}[Rounds]\label{cor:rounds}
Let $n\ge3$ and let $G$ have an edge. After any initial word, any $T\ge\lceil2n/(n-2)\rceil+2$ permutation rounds of all agents yield a word that is not universal for $G$; for $n\ge6$, five rounds suffice. Consequently, no infinite concatenation of balanced blocks of bounded size is universal at every prefix.
\end{corollary}
\begin{proof}
For rounds take $a_t=1$ in Theorem~\ref{thm:blocks}; for blocks with $a_t\le A$ the hypothesis holds as soon as $(n-2)(T-2)\ge2nA$.
\end{proof}

For round robin itself even identical rankings suffice to break EFX.

\begin{proposition}[Round robin]\label{prop:rr}
For $n\ge2$, the prefix of length $m$ of $(1\,2\cdots n)^\infty$ is universal for $K_n$ if and only if $m\le n$ or $(n,m)=(2,4)$.
\end{proposition}
\begin{proof}
For $m\le n$ all bundles are singletons, and $\A\B\A\B$ satisfies~\eqref{eq:two}. Otherwise write $m=nk+r$ with $0\le r<n$ and let all agents rank the goods chronologically. If $r=0$ (so $k\ge2$), agent $n$ sees the bundle of agent $1$ without its last good worth $\sum_{t=1}^{k-1}(m-n(t-1))$, while her own bundle is worth $\sum_{t=1}^k(m-nt+1)$; the difference is $(n-1)k-n$, positive unless $n=k=2$. If $r>0$, agent $1$ has $k+1$ goods, agent $n$ has $k$, and the difference is $(n-1)k>0$.
\end{proof}

\section{Open problems}\label{sec:open}

\begin{problem}[Sparse graphs]
Determine $\max\sigma(G)$ asymptotically for paths, cycles and stars. By Proposition~\ref{prop:snake} and Theorem~\ref{thm:linear} it is between $2n$ and $O(n)$ in all three cases. The smallest case is already intriguing: computations show that $\{0,1,\dots,18\}\subseteq\sigma(P_3)$, witnessed for instance at $m=18$ by
\[
1\,1\,2\,2\,3\,3\,2\,3\,1\,2\,2\,3\,1\,3\,3\,3\,1\,1
\]
on the path $1$--$2$--$3$, with defects $0,-1,-1,0$ on the arcs $1\to2$, $2\to1$, $2\to3$, $3\to2$, while Theorem~\ref{thm:linear} gives $\max\sigma(P_3)\le356$. What is $\sigma(P_3)$? Is it an interval?
\end{problem}

\begin{problem}[Maximal drafts]
Prove Conjecture~\ref{conj:2n}: are the snake-like words of Proposition~\ref{prop:many} the only universal words of length $2n$ for $K_n$?
\end{problem}

\begin{problem}[Drafts that may look at the preferences]
Does every positive-Borda profile admit an EFX allocation that can be obtained by sincere picking along \emph{some} order, chosen after the preferences are known? Such allocations are called sequenceable, and every Pareto-optimal allocation is sequenceable~\cite{BramsKing2005,BeynierEtAl2019Sequenceability}. So a positive answer would follow from the existence of allocations that are both EFX and Pareto optimal; but these need not exist even for positive additive valuations~\cite{MahabadiRamakrishnanTarnawski2026}. Do they exist for Borda profiles? A positive answer to either question would in particular establish the existence of EFX allocations for Borda profiles.
\end{problem}

\begin{problem}[Weaker guarantees]
If $h$ goods may be removed instead of one (EF$h$X~\cite{AkramiRezvanSeddighin2022}), how large must $h$ be, as a function of $m$ and $n$, for universal words of length $m$ to exist on connected graphs?
\end{problem}

\noindent\textbf{AI disclosure.} The authors used generative AI to help organize the manuscript, edit its exposition and check how the mathematical arguments depend on one another. The authors are responsible for the statements, proofs and references in the submitted article.

\bibliographystyle{amsplain}
\bibliography{references}

\end{document}